\documentclass[letter,12pt]{article}
\usepackage{amsfonts}
\usepackage{amssymb}
\usepackage{amsmath}
\usepackage{cite}
\newtheorem{theo}{Theorem}
\newtheorem{coro}[theo]{Corollary}

\newtheorem{lema}[theo]{Lemma}

\newcommand{\C}{\mathbb{C}}
\newcommand{\N}{\mathbb{N}}

\newcommand{\esp}{\hspace{0.1cm}}

\title{A Liv\v sic type theorem for germs of analytic diffeomorphisms}

\author{Andr\'es Navas \esp\esp \& \esp\esp Mario Ponce}
\date{}

\begin{document}

\maketitle

\vspace{-0.4cm}

\noindent{\bf Abstract.} We deal with the problem of the validity of Liv\v sic's 
theorem for cocycles of diffeomorphisms satisfying the orbit periodic obstruction 
over an hyperbolic dynamics. We give a result in the positive direction for 
cocycles of germs of analytic diffeomorphisms at the origin. 

\section{Introduction}

Given a map (dynamical system) $T \!: X \to X$ over a compact metric space $X$ and a (topological)
group $\mathcal{G}$, we consider a continuous $\mathcal{G}$-valued cocycle  $A \!: \N \times X \to
\mathcal{G}$, that is, a continuous map taking values in $\mathcal{G}$ satisfying the cocycle relation
\[
A(n+m, x)=A(n, T^mx)A(m, x)
\]
for every $m,n$ in $\N$ and every $x \in X$. This cocycle is completely determined by the continuous
function $A(\cdot) := A(1, \cdot): X\to \mathcal{G}$, and the cocycle relation yields
\[
A(n ,x)=A(T^{n-1}x)A(T^{n-2}x)\cdots A(x)
\]
for every $n \geq 1$. A natural problem is to determine conditions ensuring for such a cocycle
to be conjugated to a cocycle taking values in a ``small'' subgroup of $\mathcal{G}$. For the case
of the trivial subgroup, this property means that there exists a continuous function
$B \!: X \to \mathcal{G}$ such that
\begin{equation}\label{solv}
A(x)=B(Tx)B(x)^{-1} \quad \textrm{for all} \quad x\in X.
\end{equation}
Whenever this {\it cohomological equation} associated to the cocycle $A$ has a solution $B$, we
say that $A$ is a {\it coboundary}. The simplest obstruction for the existence of $B$ is the
{\it P(eriodic) O(orbit) O(bstruction)}: if $p\in X$ and $n\in \N$ satisfy $T^n p = p$, then
$$A(n, p) = \prod_{i=0}^{n-1} A(T^ix)
= \prod_{i=0}^{n-1} B(T^{i+1} x) B(T^i x)^{-1}
= B(T^n p)B(p)^{-1}=e_{_{\mathcal{G}}}.$$
The {\it Liv\v sic problem} consists in determining whether the 
POO($A$) condition is not only necessary
but also sufficient for $A$ being a coboundary. This terminology comes from the seminal work of
Liv\v sic \cite{LIV72}, who proved that this is the case whenever $\mathcal{G}$ is Abelian, $A$
is H\"older-continuous and $T$ is a topologically transitive hyperbolic diffeomorphism. Since then,
many extensions of this classical result have been proposed. Perhaps the most relevant is Kalinin's
recent version for $\mathcal{G} = \mathrm{GL} (d,\mathbb{C})$. In this Note, we address the
Liv\v sic problem for H\"older-continuous cocycles taking values in the group of germs of 
analytic diffeomorphisms.  In the context of general diffeomorphisms, a positive answer 
to the Liv\v sic problem is
unclear, despite several results pointing in this direction whenever a certain
localization property is satisfied. (See, for example, \cite{delallave2010}.)\\

To state our result, we denote by $\mathcal{G}erm_d$ the group of germs of local
bi-holomorphisms of the complex space $\C^d$ fixing the origin. This may be identified 
to the group of holomorphic maps $F(Z) = A_1Z+A_2Z^2+\dots$ having positive convergence 
radious, with $A_1 \in \mathrm{GL} (d, \C)$ (see \S \ref{section_germ} for the details).

\vspace{0.5cm}

\noindent{\bf Main Theorem.}
{\em Let $T: X\to X$ be a topologically transitive  
homeomorphism of a compact metric space $X$ satisfying the closing property 
(see \S \ref{remind} for the details). Let $F \!: X\to \mathcal{G}erm_d$ 
be a H\"older-continuous function/cocycle (see \S \ref{section_germ} 
for a discussion on continuity issues). If $F$ satisfies the
POO condition, then there exists a H\"older-continuous function 
$H: X \to \mathcal{G}erm_d$ such that for all $x \in X$,}
\begin{equation}\label{ecuacion_del_teorema}
F(x)=H(Tx)\circ H(x)^{-1}.
\end{equation}

\vspace{0.25cm}

This theorem should be compared with \cite{PONC11-flores}, 
where the second-named author shows a KAM-type result for 
$\mathcal{G}erm_d$-valued cocycles over a minimal torus translation.


\subsection{A remind on Liv\v sic's theorem for complex valued cocycles}
\label{remind}

Let $X$ be a compact metric space with normalized diameter 
({\em i.e.}, $diam(X) = 1$). We say that a function $f \!: X\to \C$ is
$(C, \alpha)-$H\"older-continuous for $C>0$ and $\alpha\in (0,1]$ if
for every pair of points $x, y$ in $X$,
\begin{equation}\label{holder}
|f(x)-f(y)|\leq C \esp dist_X(x, y)^{\alpha}.
\end{equation}
In the sequel, we will denote by $[f]_{\alpha}$ the smallest constant $C$ for which
$f$ is $(C, \alpha)-$H\"older-continuous. The next two results are straightforward.

\begin{lema}\label{lema_dos}
If $f$ vanishes at some point of $X$, 
then $\|f\|:=\sup_{x\in X}|f(x)|\leq [f]_{\alpha}.\quad_{\blacksquare}$
\end{lema}
\begin{lema}\label{holder_opera}
Let $f, g:X \to \C$ be two $\alpha$-H\"older-continuous functions. 
Then the functions $f+g$ and $fg$ are $\alpha-$H\"older-continuous, and
\begin{enumerate}
\item
$[f+g]_{\alpha}\leq [f]_{\alpha}+[g]_{\alpha}$.
\item $[fg]_{\alpha}\leq [f]_{\alpha}\|g\| + [g]_{\alpha} \|f\|. \quad_\blacksquare$
\end{enumerate}
\end{lema}

Let $T\!:X\to X$ be a homeomorphism and let $x, y$ be points of $X$. We 
say that the orbit segments \esp $x, Tx, \dots, T^kx$ \esp and \esp 
$y, Ty, \dots, T^ky$ \esp are {\it exponentially $\delta$-close 
with exponent \esp $\lambda>0$} \esp if for every $j=0, \dots, k$,
\[
dist_X(T^jx,T^jy)\leq \delta e^{-\lambda\min\{j, k-j\}}.
\]
We say that $T$ satisfies the {\it closing property} if there exist 
$c, \lambda, \delta_0>0$ such that for every $x\in X$ and $k\in \N$ 
so that $dist_X(x, T^kx)<\delta_0$, there exists a point $p\in X$ with 
$T^kp=p$ so that letting $\delta := c \esp dist_X(x, T^kx)$, the orbit segments 
$x, Tx, \dots, T^kx$ and $p, Tp, \dots, T^kp$ are exponentially $\delta$-close 
with exponent $\lambda$ and there exists a point $y \in X$ such that for every $j=0,\dots, k$,
\[
dist_X(T^jp, T^jy)\leq \delta e^{-\lambda j} \quad \textrm{and} 
\quad dist_X(T^jy, T^jx)\leq \delta e^{-\lambda(n-j)}.
\]
Important examples of maps satisfying the closing property are 
hyperbolic diffeomorphisms of compact manifolds.

In this work, we will use two versions of the Liv\v sic result. 
The first of these ({\em c.f.,} Theorem \ref{livsic}) is the original 
and seminal Liv\v sic theorem for complex valued cocycles. This theorem will be used in an 
iterative scheme for which having estimates for the solutions of cohomological equations will 
be relevant (see Corollary \ref{coro5}). For this reason, we review the proof and we record 
certain crucial estimates. The second version (extension) of the Liv\v sic result we will use 
({\em c.f.}, Theorem \ref{kal}) corresponds to a recent and remarkable theorem by B. Kalinin, 
who proves the Liv\v sic theorem for matrix-valued cocycles (satisfying no localization condition).

\begin{theo}[Liv\v sic, see \cite{LIV72}]\label{livsic}
Let $T \!: X\to X$ be a topologically transitive homeomorphism of a compact metric space $X$ 
satisfying the closing property. Let $\psi\!: X\to \C$ be an $\alpha-$H\"older-continuous 
function for which the POO holds, that is, for every point $p \in X$ and $k \geq 1$ such 
that $T^k p = p$, one has \esp $\sum_{j=0}^{k-1}\psi(T^{j}p)=0$. \esp Then there exists 
an $\alpha-$H\"older-continuous function $\phi:X\to \C$ that is a solution to the 
cohomological equation
\[
\phi\circ T-\phi =\psi.
\]
\end{theo}

\noindent{\it Proof.} Let $x_0\in X$ be such that $\overline{\{T^n x_0\}_{n\in \N}}=X$. We 
define $\phi$ by letting $\phi(x_0) := 0$ and $\phi(T^n x_0) := \sum_{j=0}^{n-1}\psi(T^j x_0)$. 
We next check that $\phi$ is $\alpha-$H\"older-continuous on $\{T^n x_0\}_{n\in \N}$.  
Let $n>m$. There are two cases to consider:
\begin{itemize}
\item
Assume that \esp $dist_X(T^mx_0, T^nx_0)<\delta_0$. \esp Then there exists a point 
$p\in X$ satisfying $T^{n-m}p=p$ and such that for every $j=0,\dots, n-m$,
\[
dist_X(T^j(T^mx_0), T^jp)\leq c \esp dist_X(T^nx_0, T^mx_0)e^{-\lambda\min \{j, n-m-j\}}.
\]
This yields

\begin{eqnarray*}
|\phi(T^nx_0)-\phi(T^mx_0)|&=&\left|\sum_{j=0}^{n-m-1}\psi(T^{m+j}x_0)\right|\\
&=&\left|\sum_{j=0}^{n-m-1}\left(\psi(T^{m+j}x_0)-\psi(T^jp)\right)+\sum_{j=0}^{n-m-1}\psi(T^jp)\right|\\
&\leq&\sum_{j=0}^{n-m-1}\left|\psi(T^{m+j}x_0)-\psi(T^jp)\right|\\
&\leq&\sum_{j=0}^{n-m-1}[\psi]_{\alpha} \esp dist_X(T^{m+j}x_0, T^jp)^{\alpha}\\
&\leq&\sum_{j=0}^{n-m-1}c^{\alpha}[\psi]_{\alpha} \esp 
dist_X(T^nx_0, T^mx_0)^{\alpha} \esp e^{-\lambda\alpha\min \{j, n-m-j\}}\\
&\leq&\frac{2 \esp c^{\alpha} [\psi]_{\alpha}}{1-e^{-\lambda\alpha}} \esp 
dist_X(T^nx_0, T^mx_0)^{\alpha}.
\end{eqnarray*}
\item Assume that $dist_X(T^nx_0, T^mx_0)\geq  \delta_0$. Since $x_0$ has dense orbit and 
$X$ is compact, there exists $N\in \N$, depending only on $X, T$, and $\delta_0$, such 
that $\{x_0, Tx_0, \dots, T^Nx_0\}$ is a $\delta_0$-dense set in $X$. For $n-m\leq N$, 
one easily shows that
\[
|\phi(T^nx_0)-\phi(T^mx_0)|\leq N\|\psi\|.
\]
For $n-m>N$, there exist $r, s$ in $\{0,1,\dots, N\}$ 
such that $dist_X(T^{s}x_0, T^nx_0)\leq \delta_0$ and 
$dist_X(T^{r}x_0, T^mx_0)\leq \delta_0$. Using the preceding case, 
this yields 
\begin{eqnarray*}
|\phi(T^nx_0)-\phi(T^mx_0)|
\!\!&\leq&\!\!
|\phi(T^nx_0)-\phi(T^{s}x_0)|+|\phi(T^mx_0)-\phi(T^{r}x_0)|
+|\phi(T^{s}x_0)-\phi(T^{r}x_0)|\\
\!\!&\leq&\!\! 
\frac{4[\psi]_{\alpha}c^{\alpha}}{1-e^{-\lambda\alpha}}\delta_0^{\alpha}+N\|\psi\|\\
\!\!&\leq&\!\! 
\left(\frac{4[\psi]_{\alpha}c^{\alpha}}{1-e^{-\lambda\alpha}}
+\frac{N\|\psi\|}{\delta_0^{\alpha}}\right)dist_{X}(T^nx_0, T^mx_0)^{\alpha}.\quad_{\blacksquare}
\end{eqnarray*}
\end{itemize}

A careful reading of the proof above yields useful estimates enclosed in the next 

\begin{coro}\label{coro5}
The solution $\phi$ to the cohomological equation is $\alpha$-H\"older 
continuous, and there exists $K$ depending only on $T, X,$ and $\alpha$ such 
that $[\phi]_{\alpha}\leq K([\psi]_{\alpha}+\|\psi\|). \quad_{\blacksquare}$
\end{coro}

\begin{theo}[Kalinin, see \cite{KAL10}]\label{kal}
Let $T$ be a topologically transitive homeomorphism of a compact metric space $X$ satisfying 
the closing property. Let $A:X\to GL(d, \C)$ be an $\alpha$-H\"older function for which the 
$POO(A)$ holds. Then there exists an $\alpha$-H\"older function $C:X\to GL(d, \C)$ such that
for all $x \in X$,
\[
A(x) = B(Tx) B(x)^{-1}.
\]
\end{theo}


\subsection{The group $\mathcal{G}erm_d$} 
\label{section_germ}

For $d \geq 1$, we introduce the following (classical) notation:
\begin{itemize}
\item${\bf j}:=(j_1, \dots, j_d)$ is a positive integer lattice point, 
with $j_i\geq 0$ for every $1\leq i\leq d$. 
\item$|{\bf j}|:=j_1+\dots+j_d$.
\item${\bf j} \preceq {\bf k}$ \esp if \esp $j_i\leq k_i$ for every $1\leq i\leq d$. 
\item${\bf j} \prec {\bf k}$ \esp if \esp ${\bf j} \preceq {\bf k}$ 
and $j_{i_*}<k_{i_*}$ for some $i_*$.
\item$Z=(z_1, z_2, \dots, z_d)$ is a point in $\C^d$.
\item$Z^{\bf j}:=z_1^{j_1}z_2^{j_2}\cdots z_d^{j_d}$.
\end{itemize}
Then we can define a formal power series on $\C^d$ as 
$
F(Z) := (F_1(Z), F_2(Z), \dots, F_d(Z)),
$
where each $F_i(Z)$ has the form
\[
F_i(Z)=\sum_{{\bf j}\geq 0} t_{\bf j}^iZ^{\bf j}
\]
for some coefficients $t_{\bf j}^i\in \C$. This formal power series becomes an analytic map if 
there exists $R>0$ such that \esp $\limsup_{\bf j}|t_{\bf j}^i|^{\frac{1}{|{\bf j}|}} \leq 
\frac{1}{R}$ \esp for every $i$. Indeed, in this case, each $F_i$ is a convergent series on  
$D(0, R)^d$ (that is, for $Z \!=\! (z_1,\ldots,z_d)$ such that $|z_s| \!<\! R$ 
holds for every $s$).

Let $\mathcal{H}(d, R)$ be the set of continuous functions 
$F \!: \overline{D(0, R)^d }\to \C^d$ that are convergent power 
series in $D(0, R)^d$ and satisfy $F'(0) \!\in\! GL(d , \C)$. 
We endow this complex vector space with the inner product
\[
\langle F, G \rangle_R 
:= \sum_i\left(\int_{\partial D(0,R)^d}F_i\overline{G_i} \esp dZ\right).
\]
The $L^2$-norm of an element $F \! \in \! \mathcal{H}(d, R)$ of the form 
$F_i(Z)=\sum_{|{\bf j}|\geq 0}t^i_{\bf j}Z^{\bf j}$ is 
\[
\|F\|_{2, R} := \langle F, F \rangle_R^{1/2} 
= \left(\sum_i\sum_{|{\bf j}|\geq 1}|t_{\bf j}^i|^2R^{2|{\bf j}|}\right)^{1/2}.
\]
We let $\mathcal{H}_0 (d, R)$ be the subset of $\mathcal{H} (d,R)$ formed by those $F$ 
satisfying $F(0) = 0$, and we let the set of {\it local holomorphic diffeomorphisms} 
of $\C^d$ be defined as
\[
\mathcal{G}_d := \bigcup_{R> 0}\mathcal{H}_0(d, R).
\]
On this set, we introduce the following equivalence relation: We say that 
$F, G$ in $\mathcal{G}_d$ are equivalent if there exists a neighborhood 
of the origin on which $F$ and $G$ coincide. Under this identification, 
the set $\mathcal{G}_d$ becomes a group, that we call the {\it group 
of germs of analytic diffeomorphisms} of $\C^d$ and we denote by 
$\mathcal{G}erm_d$.\\

Although we will not worry about giving a topology on $\mathcal{G}erm_d$, we will 
certainly need to consider maps from $X$ to $\mathcal{G}erm_d$ that are ``continuous'' 
in some precise sense. Since $X$ is compact, any reasonable definition should lead 
to functions that factor throughout an space $\mathcal{H}_0(d, R)$ for some 
positive $R$. Accordingly, given $C>0$, $\alpha \in (0,1]$, and $R>0$, 
a map $\Psi \!: X \to \mathcal{H}_0(d, R)$ will be said to be 
$\left(C, \alpha, R\right)-$H\"older-continuous if 
$\Psi(x)$ belongs to $\mathcal{H}_0(d, R)$ for every 
$x \in X$, and for every pair of points $x, y$ in $X$,
\[
\|\Psi(x)-\Psi(y)\|_{2, R} \leq C \esp dist_X(x, y)^{\alpha}.
\]
In terms of the coefficients of the power series, this condition reads as follows:

\vspace{0.1cm}

\begin{lema}\label{lema6666}
If $\Psi \!: X \to \mathcal{H}_0(d, R)$ is $\left(C, \alpha, R\right)-$H\"older 
and writes as 
$$\Psi_i(x)(Z)=\sum_{|{\bf j}|>0} t^i_{\bf j}(x)Z^{\bf j},$$ 
then each coefficient $t^i_{\bf j}:X\to \C$ is a 
$\left(\frac{C}{R^{|{\bf j}|}}, \alpha\right)$-H\"older-continuous function.
\end{lema}

\noindent{\it Proof.} The H\"older condition for $\Psi$  yields
\[
\left(\sum_i\sum_{|{\bf j}|\geq 1}|t_{\bf j}^i(x)-t^i_{\bf j}(y)|^2R^{2|{\bf j}|}\right)^{1/2} 
\leq C \esp dist_X(x, y)^{\alpha},
\]
which implies that
\[
|t_{\bf j}^i(x)-t^i_{\bf j}(y)|^2 \leq \frac{C^2}{R^{2|{\bf j}|}} 
dist_X(x, y)^{2\alpha}.\quad_{\blacksquare}
\]

\vspace{0.2cm}

In an opposite direction, given a list 
$\{t_{\bf j}^i:X\to \C, \esp {\bf j} \succeq 0, \esp 1\leq i\leq d \}$ 
of continuous functions, we are interested in finding conditions 
ensuring that $F := (F_1, \dots, F_d)$ formally defined by 
$F_i(x)(Z):=\sum_{{\bf j}}t^i_{\bf j}(x)Z^{\bf j}$ represents a 
convergent power series lying in $\mathcal{H}_0(d, R)$ for some $R>0$. 

\vspace{0.1cm}

\begin{lema}\label{lema777}
Assume that each function $t_{\bf j}^i$ is a 
$\left(\frac{C}{R^{|{\bf j}}|}, \alpha\right)$-H\"older-continuous 
function for some positive constants $C,R$. Assume also that each $t_{\bf j}^i$ 
vanishes at some point of $X$. Then for all $\delta \!<\! 1$, the formal power series $F_i$ 
above is convergent on $D(0, R)^d$, and $x\mapsto F(x) \! = \! (F_1(x), \dots, F_d(x))$ is a 
$\left(O\!\left(\frac{\delta}{1-\delta}\right)^{1/2}\!, \esp \alpha\right)$-H\"older 
continuous map from $X$ to $\mathcal{H}_0(d, \delta R)$.
\end{lema}

\noindent{\it Proof.} Since each $t_{\bf j}^i$ vanishes at some point of $X$, Lemma \ref{lema_dos} 
gives $\|t_{\bf j}^i\|\leq \frac{C}{R^{|{\bf j}|}}$ for every $i, {\bf j}$. This implies that 
each $F_i$ is a convergent power series on $D(0, R)^d$. Moreover, for all $x, y$ in $X$,

\begin{eqnarray*}
\|F(x)-F(y)\|_{2, \delta R}^2 
&=& \sum_i \sum_{\bf j}|t^i_{\bf j}(x)-t^i_{\bf j}(y)|^2(\delta R)^{2|{\bf j}|}\\
&\leq& \sum_i\sum_{\bf j}C^2dist_X(x, y)^{2\alpha}\delta ^{2|{\bf j}|}\\
&=& d \esp C^2dist_X(x, y)^{2\alpha}\sum_{s=1}^{\infty}\sum_{ |{\bf j}|=s} \delta^{2s}\\
&=& d \esp C^2dist_X(x, y)^{2\alpha}\sum_{s=1}^{\infty}\frac{(s+d-1)!}{s!(d-1)!}\delta^{2s}\\
&=& d \esp C^2 \esp O\!\left(\frac{\delta}{1-\delta}\right)dist_X(x, y)^{2\alpha}.
\quad_{\blacksquare}
\end{eqnarray*}

\paragraph{The Faa di Bruno formula.} We will need to consider compositions of power series in 
several complex variables. The following is a  simplified formulation of the multivariate 
version by Constantine and Savits \cite{CONSSAVI96} of the famous Faa di Bruno formula:

\begin{theo}[see \cite{CONSSAVI96}]
Let $A(Z)=\sum_{|{\bf j}|\geq 1}a_{\bf j}Z^{\bf j}$ and 
$B_i(Z)=\sum_{|{\bf j}|\geq 1}b^i_{\bf j}Z^j, \esp 1 \leq i \leq d$,
be formal power series in $d$ variables. Then the power series 
\[
C(Z)=A\left(B_1(Z), B_2(Z), \dots, B_d(Z)\right)=\sum_{|{\bf j}|\geq 1}c_{\bf j}Z^{\bf j}
\]
has coefficients
\begin{equation}\label{faadi}
c_{\bf j_*}=\sum_{|\bf j|=1}a_{\bf j}b_{\bf j_*}^{\bf j} \esp 
+ \! \sum_{1<|{\bf j}|, \ {\bf j}\leq {\bf j_*}}a_{\bf j}P({\bf j_*}, {\bf j})\{B\}, 
\end{equation}
where $P({\bf j_*}, {\bf j})\{B\}$ is  polynomial in the variables 
$\{b_{\bf \tilde{j}}^i\}^{1\leq i\leq d}_{{\bf \tilde j}<{\bf j_*}}$ 
that is homogeneous of degree $|{\bf j}|$ and has positive integer coefficients.
\end{theo} 

The Faa di Bruno formula is actually much more precise but requires hard notation. 
For instance, in the case $d=1$, one has
\[
P(j_*, j)\{B\}=\sum_{r_1+\dots+r_j=j_*}B_{r_1}\cdots B_{r_j}.
\]


\paragraph{A generating function.} Let us define 
$J \!:\! D(0,1)^d\to \C^d$ by the convergent power series
\[
J_i(Z)=z_i-\sum_{|{\bf j}|> 1}Z^{\bf j}.
\]
Since $DJ(0)=id_{\C^d}$, there exists an analytic map $G$ defined 
in a neighborhood of the origin in $\C^d$ such that $G(0) = 0$ and 
\begin{equation}\label{misma_formal}
J \circ G (Z) = Z \esp \esp \mbox{ for every } Z \mbox{ in that neighborhood }. 
\end{equation}
In terms of power series, one can write
\[
G_i(Z)=z_i+\sum_{|{\bf j}|>1} g^i_{\bf j}Z^{\bf j},
\]
where the coefficients verify $|g^i_{\bf j}|<K^{|{\bf j}|-1}$ for some $K>0$ and  
every $|{\bf j}| \succ 1$. Moreover, these coefficients satisfy a fundamental recurrence 
relation. Indeed, using $J\circ G(Z) = Z$ and the Faa di Bruno formula (\ref{faadi}),  
one obtains 
\begin{equation}\label{recurrencia}
0 \esp = \esp g^i_{\bf j_*} \esp - 
\sum_{1<|{\bf j}|, \ {\bf j}\leq {\bf j_*}}P({\bf j_*}, {\bf j})\{G\}.
\end{equation}
Recall that $P({\bf j_*}, {\bf j})\{G\}$ depends only on the values of 
$g^s_{\bf \tilde j}$ for ${\bf \tilde j} \prec {\bf j_*}$ and every $s$. 
Hence, one can recursively compute $g_{\bf j_*}^i$ in terms of the 
previously defined $g^s_{\bf \tilde j}$.

For any $S>0$, we consider 
$J_S:D(0, S^{-1})^d\to \C^d$ defined by $J_S(Z):=\frac{1}{S}J(SZ)$. 
When solving the equation \esp $J_S\circ G_S(Z)=Z$, \esp one gets a  map 
$G_S = (G_{S, 1}, \dots, G_{S, d})$, where each $G_{S, i}(Z)$ has the form 
\esp $G_{S, i} (Z) = z_i + \sum_{|{\bf j}|>1} g_{S, {\bf j}}^iZ^{\bf j}$ 
\esp for certain coefficients $g_{S, {\bf j}}^i$ satisfying
\begin{equation}\label{recurrenciaS}
g^i_{S, {\bf j_*}}=\sum_{1<|{\bf j}|, \ {\bf j} 
\leq {\bf j_*}}S^{|\bf j|-1}P({\bf j_*}, {\bf j})\{G_S\}.
\end{equation}

\begin{lema}\label{lema_de_S}
Each coefficient $g^i_{S, {\bf j}}$ is a positive real number. 
Moreover, there exists a constant $\mathcal{R}=\mathcal{R}(S)>0$ such that 
$g^i_{S, {\bf j}}\leq \mathcal{R}^{|{\bf j}|-1}$ for every 
${\bf j}$. $\quad_{\blacksquare}$
\end{lema}


\subsection{Proof of the Main Theorem}

\paragraph{A first reduction.} 
Let \esp $F(x)(Z)=A_1(x)Z + 
\left(\sum_{|{\bf j}|>1} a^i_{\bf j}(x)Z^{\bf j}\right)_{1\leq i\leq d}$ 
be the power series expansion of the cocycle viewed as a \esp 
$\left(C, \alpha, R\right)-$H\"older-continuous function 
$\Psi \!: X \to \mathcal{H}_0(d, R)$. The map 
$x\mapsto A_1(x)\in GL(d, \C)$ is an $\alpha$-H\"older-continuous 
function. Since POO($F$) holds, we must have
\[
\prod_{j=0}^{n-1}A_1(T^jp) \esp\esp
= \esp\esp \frac{\partial}{\partial Z}F(T^{n-1}p)\circ \dots \circ F(p)\bigg|_{Z=0} 
= \esp\esp id_{\C^d}
\]
for every $p\in X$ and $n\in \N$ such that $T^np=p$. In other words, the $GL(d, \C)$-valued 
cocycle $A_1$ satisfies the POO. By Kalinin's version of the Liv\v sic theorem, there exists 
an $\alpha$-H\"older-continuous function $H_1:X\to GL(d, \C)$ such that $A_1(x)=H_1(Tx)H_1(x)^{-1}$ 
for all $x\in X$. Consequently, the $\mathcal{G}erm_d$-valued cocycle \esp $H_1(x)(Z):=H_1(x)Z$ 
\esp conjugates $F$ to a cocycle of the form 
$$(x, Z)\longmapsto 
\left(Tx, Z+\left(\sum_{|{\bf j}|>1}a^i_{\bf j}(x)Z^{\bf j}\right)_{1\leq i\leq d}\right).$$ 
Thus, we can assume that $A_1(x)=id_{\C^d}$ for all $x\in X$.

\paragraph{An iterative procedure.} We look for a map $H \!: X\to \mathcal{G}erm_d$ 
solving the cohomological equation (\ref{ecuacion_del_teorema}) having the form 
$H(x)(Z)=Z+\left(\sum_{|{\bf j}|>1}h^i_{\bf j}(x)Z^{\bf j}\right)_{1\leq i\leq d}$. 
Notice that this equation may be written as \esp $F(x)\circ H(x)=H(Tx)$. \esp 
Applying the Faa di Bruno formula (\ref{faadi}) to the left-side expression, 
one concludes that each coefficient $h_{\bf j}^{i}$ can be defined recursively 
as the solution of a cohomological equation for a $\C$-valued data:
\begin{displaymath}
\begin{array}{crcl}
(ec_{\bf j_*}^{i}) & \quad h^{i}_{\bf j_*}(Tx)-h^{i}_{\bf j_*}(x) 
&=& 
\sum_{1<|{\bf j}|, \ {\bf j}\leq {\bf j_*}}a^{i}_{\bf j}(x)P({\bf j_*}, {\bf j})\{H\}(x).
\end{array}
\end{displaymath}
A necessary  condition for the existence of the coefficient 
$h^{i}_{\bf j_*}$ is that the POO condition holds for the function
\begin{equation}\label{formularecursiva}
R^i_{\bf j_*} \esp
:= \sum_{1<|{\bf j}|, \ {\bf j}\leq {\bf j_*}}a^{i}_{\bf j}P({\bf j_*}, {\bf j})\{H\}.
\end{equation}

\vspace{0.1cm}

\begin{lema}
Each $R_{\bf j_*}^i$, with $i, |{\bf j_*}| \succ 1$, is a well-defined 
$\alpha$-H\"older-continuous function for which the POO holds. As a consequence, 
given any $x_0 \in X$, the equation $(ec_{\bf j_*}^{i})$ has an $\alpha$-H\"older-continuous 
solution $h^i_{\bf j_*}$ vanishing at $x_0$.
\end{lema}

\noindent{\it Proof.} Suppose that the conclusion of the lemma holds for every 
${\bf j}$ such that $|{\bf j}|<k$, and let us consider the case where ${\bf j} = k$. 
Using the explicit formula (\ref{formularecursiva}), Lemma \ref{holder_opera} shows 
that the function $R_{\bf j_*}^i$ is $\alpha$-H\"older-continuous. Consider the 
continuous $\mathcal{G}erm_d$-valued function
\[
H_{<k}:x\mapsto Z+\left(\sum_{|{\bf j}|<k}h^i_{\bf j}(x)Z^{\bf j}\right)_{1\leq i\leq d}.
\]
An easy computation shows that \esp 
$\tilde F(x) := H_{<k}(Tx)\circ F(x)\circ H_{<k}(x)^{-1}$ \esp has the form
\[
\tilde F(x)(Z) 
= Z+\left(\sum_{|{\bf j}|=k}R_{\bf j}^i(x)Z^{\bf j} 
+ \sum_{{|\bf j}|>k}\tilde a_{\bf j}^i(x)Z^{\bf j}\right)_{1\leq i\leq d}
\]
for some H\"older-continuous functions $\tilde a_{\bf j}^i:X\to \C$. 
Moreover, for any $x\in X$ and $m\in \N$, one has
\[
\tilde F(T^{m-1}x)\circ \dots \circ \tilde F(x)(Z) 
= Z+\left(\sum_{|{\bf j}|=k}\left(\sum_{v=0}^{m-1}R_{\bf j}^i(T^vx)\right)Z^{\bf j} 
+ \mathcal{O}(|Z|^{k+1})\right)_{1\leq i\leq d}.
\]
Since $\tilde F$ is conjugated to $F$, the POO($\tilde F$) holds. By the 
previous equality, this implies that for all $p\in X$ and $n\in \N$ such 
that $T^np=p$, one has \esp $\sum_{v=0}^{n-1} R_{\bf j}^i(T^vx)=0$. \esp 
Therefore, the POO($R_{\bf j}^i$) holds, which allows applying Livsic's 
theorem to ensure the existence of an $\alpha$-H\"older-continuous 
solution to $(ec_{\bf j_*}^{i})$. Finally, by adding a constant if 
necessary, we may assume that this solution vanishes at $x_0$.  
$\quad_{\blacksquare}$\\
 
To prove that the (up to now) formal map $H$ is a genuine local diffeomorphism 
(that is, each formal power series $Z \mapsto z_i+\sum_{|{\bf j}|>1}h_{\bf j}^i(x)Z^{\bf j}$ 
is convergent in a certain (uniform) neighborhood of the origin), we will need to estimate 
the growth of the $\alpha$-H\"older constant of the coefficients $h_{\bf j}^i$. Indeed, if 
we show that this growth is at most exponential, then Lemma \ref{lema777} will apply, thus 
concluding the proof of the Main Theorem. To get the desired control, we will use the {\it majorant 
series method} introduced by Siegel in his treatement \cite{SIEG42} of the linearization 
theorem for holomorphic germs with Diophantine rotation number (see also \cite{STEN61} for the 
higher-dimensional case).  

\begin{lema}\label{tamanodeh}
There exists $S>0$ such that
\[
[h_{\bf j}^i]_{\alpha}\leq g_{S, {\bf j}}^i
\]
for every ${\bf j}, i$, where $h_{S, {\bf j}}^i$ is defined as in (\ref{recurrenciaS}). 
Consequently, $\|h_{\bf j}^i\|$  grows at most exponentially.
\end{lema}

\noindent {\it Proof.} Since $F$ takes values on some $\mathcal{H}_0(d, R)$ and is a 
$\alpha$-H\"older function, there exists $\kappa>0$ such that
\[
\|a^i_{\bf j}\|\leq \kappa^{|{\bf j}|}\quad \textrm{ and} \quad 
\ [a^i_{\bf j}]_{\alpha}\leq \kappa^{|{\bf j}|}.
\]
Assume that $[h_{\bf j}^i]_{\alpha} \leq g_{S, {\bf j}}^i$ 
for every ${\bf j} \preceq {\bf j_*}$. Since $h_{\bf j}^i$ 
vanishes at $x_0$ (except for $|\bf j|=1$, for which $h_{\bf j}^i\equiv 1$), 
we also have $\|h_{\bf j}^i\|\leq g_{S, {\bf j}}^i$ for every 
${\bf j} \preceq {\bf j_*}$. Moreover, since 
$P({\bf j_*},{\bf j})\{H\}$ is an homogeneous polynomial in 
$\{h_{\tilde{\bf j}}^s\}^{1\leq s\leq d}_{\tilde {\bf j} < {\bf j_*}}$
with positive coefficients, 
\[
\|P({\bf j_*},{\bf j})\{H\}\|\leq P({\bf j_*},{\bf j})\{\|H\|\} 
\leq P({\bf j_*},{\bf j})\{G_S\}.
\]
Except for $|\bf j|=1$ (for which $h_{\bf j}^i\equiv 1$), 
every $h^i_{\bf j}$ vanishes at $x_0$. Therefore, by 
Lemma \ref{holder_opera},
\[
\left[P({\bf j_*},{\bf j})\{H\}\right]_{\alpha} 
\leq 2^{|\bf j|-1}P({\bf j_*},{\bf j})\{G_S\}.
\]
The fundamental estimate of Corollary \ref{coro5} then yields
\begin{eqnarray*}
[h^i_{\bf j_*}]_{\alpha}&\leq &K\left(\left[\sum_{{\bf j} 
\leq {\bf j_*}}a^i_{\bf j}P({\bf j_*}, {\bf j})\{H\}\right]_{\alpha}
+\left\|\sum_{{\bf j}\leq {\bf j_*}}a^i_{\bf j}P({\bf j_*}, {\bf j})\{H\}\right\|\right)\\
&\leq & 
K\left(\sum_{{\bf j}\leq {\bf j_*}}\|a^i_{\bf j}\|[P({\bf j_*}, {\bf j})\{H\}]_{\alpha}+\sum_{{\bf j} 
\leq {\bf j_*}}[a^i_{\bf j}]_{\alpha}\|P({\bf j_*}, {\bf j})\{H\}\|+\sum_{{\bf j}
\leq {\bf j_*}}\|a^i_{\bf j}\|\|P({\bf j_*}, {\bf j})\{H\}\|\right)\\
&\leq& 
\sum_{{\bf j}\leq {\bf j_*}}K\left((2\kappa)^{|{\bf j}|} 
+2\kappa^{|{\bf j}|}\right)P({\bf j_*}, {\bf j})\{G_S\}\\
&<& 
g^i_{S, {\bf j_*}},
\end{eqnarray*}
where the last inequality holds by taking $S \gg 2K\kappa$.$\quad_{\blacksquare}$\\

\vspace{0.25cm}


\begin{small}

\noindent{\bf Acknowledgments.} Both authors were funded by the Math-AMSUD Project DySET. 
Mario Ponce was also funded by the Fondecyt Grant 11090003.

\end{small}


\begin{footnotesize}

\vspace{0.25cm}

\noindent{Andr\'es Navas}

\noindent{Dpto de Matem\'atica y C.C., USACH}

\noindent{Alameda 3363, Estaci\'on Central, Santiago, Chile}

\noindent{E-mail: andres.navas@usach.cl}\\

\noindent{Mario Ponce}

\noindent{Facultad de Matem\'aticas, PUC}

\noindent{Casilla 306, Santiago 22, Chile}

\noindent{E-mail: mponcea@mat.puc.cl}

\end{footnotesize}
\end{document}